\documentclass[preprint,10pt,twocolumn,5p,authoryear]{elsarticle}

\usepackage{amsmath,amssymb,amsfonts}
\usepackage{dsfont}
\usepackage{mathrsfs}

\usepackage{xcolor}

\usepackage{bbm}
\usepackage{xspace}

\newcommand{\psh}[2]{\ensuremath{\langle #1,#2\rangle}\xspace}
\DeclareMathAlphabet{\pazocal}{OMS}{zplm}{m}{n}
 
\newtheorem{thm}{Theorem}[section]

\newtheorem{lemma}{Lemma}[section]

\newtheorem{prop}{Proposition}[section]

\newtheorem{defn}{Definition}[section]

\newtheorem{rem}{Remark}[section]

\newproof{pf}{\textbf{Proof}}
\renewcommand*{\qed}{\hfill\ensuremath{\blacksquare}}
\newcommand{\vertiii}[1]{{\vert\kern-0.25ex\vert\kern-0.25ex\vert #1 
    \vert\kern-0.25ex\vert\kern-0.25ex\vert}}

\journal{European Journal of Control}

\begin{document}

\begin{frontmatter}

\title{On Dirac structure of infinite-dimensional stochastic port-Hamiltonian systems}

\author[Luxembourg,Namur]{François Lamoline}\ead{francois.lamoline@uni.lu}            
\author[Namur]{Anthony Hastir}\ead{anthony.hastir@unamur.be} 

\address[Luxembourg]{University of Luxembourg, Luxembourg Centre for Systems Biomedicine, Avenue du Swing 6, L-4367 Belvaux, Luxembourg}    
\address[Namur]{University of Namur, Department of Mathematics and Namur Institute for Complex Systems (naXys), Rue de Bruxelles 61, B-5000 Namur, Belgium}



\begin{abstract}       
Stochastic infinite-dimensional port-Hamiltonian systems (SPHSs) with multiplicative Gaussian white noise are considered. In this article we extend the notion of Dirac structure for deterministic distributed parameter port-Hamiltonian systems to a stochastic ones by adding some additional stochastic ports. Using the Stratonovich formalism of the stochastic integral, the proposed extended interconnection of ports for SPHSs is proved to still form a Dirac structure. This constitutes our main contribution. We then deduce that the interconnection between (stochastic) Dirac structures is again a (stochastic) Dirac structure under some assumptions. These interconnection results are applied on a system composed of a stochastic vibrating string actuated at the boundary by a mass-spring system with external input and output. This work is motivated by the problem of boundary control of SPHSs and will serve as a foundation to the development of stabilizing methods.
\end{abstract}

\end{frontmatter}

\section{Introduction}
Linear distributed port-Hamiltonian systems constitute a powerful class of systems for the modelling, the analysis, and the control of distributed parameter systems. It enables us to model many physical systems such as beam equations, transport equations or wave equations, see for instance \cite{zwart}. A comprehensive overview of the literature on this class of systems can be found in \cite{Hashem20}. In order to cover an even larger set of systems that admits a port-Hamiltonian representation, some authors have defined the notion of dissipative or irreversible port-Hamiltonian systems, see e.g. \cite{Mora2021, Mora2021Bis, Caballeria} and \cite{Ramirez2022}. Port-Hamiltonian systems (PHSs) are characterized by a Dirac structure and an Hamiltonian. Dirac structures consist of the power-preserving interconnection of different ports elements and were first introduced in the context of port-Hamiltonian systems in \cite{vds02}. It was then extended for higher-order PHSs in \cite{Legorrec2005} and \cite{Villegas}. A fundamental property of Dirac structures is that the composition of Dirac structures still forms a Dirac structure, provided some assumptions. This induces the main aspect of the port-Hamiltonian modelling, which is that the power-conserving interconnection of PHSs is still a port-Hamiltonian system.

Stochastic models are powerful to take into account neglected random effects that may occur when working with real plants. Especially, random forcing, parameter uncertainty or even boundary noise can impact the behavior of dynamical systems. In particular, PHSs interact with their environment through external ports, which can be a cause of randomness in many different ways as explained in \cite{LamolineEjc}. The stochastic extension of port-Hamiltonian systems was first proposed in \cite{Lazaro2008} for Poisson manifolds. On finite-dimensional spaces the class of nonlinear time-varying stochastic port-Hamiltonian systems (SPHSs) was presented in \cite{satoh}. A stochastic extension of distributed port-Hamiltonian systems was first developed in \cite{Lamolinethesis} and in \cite{LamolineIEEE}. In addition the passivity property of SPHSs was investigated in \cite{Lamoline-cdc} and \cite{LamolineEjc}. However, in these works only a state space representation described by a stochastic differential equation (SDE) is given. As far as known, few efforts have been done for describing the underlying geometric structure of SPHSs. One can cite \cite{cordoni2019stochastic}, where finite-dimensional stochastic port-Hamiltonian systems are modelled using the Stratonovich and Ito formalisms.

In this paper the notion of Dirac structure with stochastic port-variables is explored. Our central idea consists in extending the original Dirac structure of deterministic first-order PHSs by adding further noise ports to the port-based structure. These specific noise ports are devoted to represent the interaction of the dynamical system with its random environment. In order to preserve the power-preserving interconnection we consider the Stratonovich formulation of the stochastic integral, see for instance \cite{duan2014}. Interested readers may also be referred to \cite{CurtainPritchard} and \cite{daprato} for further details on infinite-dimensional SDEs. 

The content of this article is as follows. In Section \ref{sec:background} we introduce the basic concepts on Dirac structures together with the class of deterministic port-Hamiltonian systems. In Section \ref{sec:SPHSs} a port-based representation for SPHSs is presented and it is shown to form a Dirac structure, which is the main contribution of the paper. Section \ref{InterSection} is dedicated to the illustration of our central result, by showing that some interconnection between the newly defined Dirac structure and another arbitrary Dirac structure that shares common ports is still a Dirac structure. A stochastic damped vibrating string actuated by a mass-spring system at the boundary is then presented as an example. We conclude and discuss some future works in \ref{sec:concl}.

\section{Background on Dirac structure}\label{sec:background}
In this section we introduce some notions on distributed port-Hamiltonian systems, Tellegen structures and Dirac structures. Let us first recall the definitions of Tellegen and Dirac structures for linear distributed PHSs, see e.g. \cite{vds02}, \cite{Legorrec2005} and \cite{Kurula}. Let $\mathcal{E}$ and $\mathcal{F}$ be two Hilbert spaces endowed with the inner products $\psh{\cdot}{\cdot}_\mathcal{E}$ and $\psh{\cdot}{\cdot}_{\mathcal{F}}$, respectively. The spaces $\mathcal{E}$ and $\mathcal{F}$ denote the effort and the flow spaces, respectively. We define the bond space $\mathcal{B}:=\mathcal{F}\times\mathcal{E}$ equipped with the following inner product
\begin{equation}
\left\langle \left(\begin{matrix}
f_1\\
e_1
\end{matrix}\right),\left(\begin{matrix}
f_2\\
e_2
\end{matrix}\right)\right\rangle_{\mathcal{B}} = \psh{f_1}{f_2}_{\mathcal{F}} + \psh{e_1}{e_2}_{\mathcal{E}}  
\end{equation}
for all $(f_1, e_1), (f_2, e_2) \in \mathcal{B}$.\\ 
To define Tellegen or Dirac structures, the bond space is endowed with the bilinear symmetric pairing given by 
\begin{equation}
\left\langle \left(\begin{matrix}
f_1\\
e_1
\end{matrix}\right),\left(\begin{matrix}
f_2\\
e_2
\end{matrix}\right)\right\rangle_{+} = \psh{f_1}{j^{-1}e_2}_\mathcal{F} + \psh{e_1}{j f_2}_\mathcal{E},
\label{SymmetricPairing}
\end{equation}   
with $j:\mathcal{F}\to\mathcal{E}$ being an invertible linear mapping. The bilinear pairing $\left\langle \cdot,\cdot\right\rangle_{+}$ represents the power.\\ 
Let $\mathcal{V}$ be a linear subspace of $\mathcal{B}$. The orthogonal subspace of $\mathcal{V}$ with respect to the bilinear pairing $\left\langle \cdot,\cdot\right\rangle_{+}$ is defined as 
\begin{equation}
\mathcal{V}^\perp:= \lbrace b\in \mathcal{B} : \psh{b}{v}_+ = 0, \text{	for all	} v\in\mathcal{V}\rbrace.
\end{equation} 
These tools enable us to define Tellegen and Dirac structures, see \cite[Definition 2.1]{Kurula}.
\begin{defn}\label{Tellegen}
A linear subspace $\mathcal{D}$ of the bond space $\mathcal{B} := \mathcal{F}\times\mathcal{E}$ is called a Tellegen structure if $\mathcal{D}\subset \mathcal{D}^\perp$, where the orthogonal complement is understood with respect to the bilinear pairing $\langle\cdot,\cdot\rangle_{+}$, see (\ref{SymmetricPairing}).
\end{defn}

\begin{defn}\label{diracdef}
A linear subspace $\mathcal{D}$ of the bond space $\mathcal{B}$ is said to be a Dirac structure if 
\begin{equation}
\mathcal{D}^\perp= \mathcal{D}.
\label{dirac_cond}
\end{equation}
\end{defn}
Note that the condition (\ref{dirac_cond}) implies that the power of any element of the Dirac structure is equal to zero, i.e., 
\begin{equation*}
\left\langle \left(\begin{matrix}
f\\
e
\end{matrix}\right),\left(\begin{matrix}
f\\
e
\end{matrix}\right)\right\rangle_{+} = 2 \psh{f}{j^{-1}e}_\mathcal{F} = 0,
\end{equation*}
for any $(f,e)\in \mathcal{D}$, where the relation $\langle f,j^{-1}e\rangle_\mathcal{F} = \langle jf,e\rangle_\mathcal{E}$ has been used. The underlying structure of port-Hamiltonian systems forms a Dirac structure, which links the port-variables in a way that the total power is equal to zero.
A distributed port-Hamiltonian system is described by the following partial differential equation 
\begin{equation}
\begin{split}
\dfrac{\partial \varepsilon}{\partial t}(\zeta ,t) = P_1 \dfrac{\partial}{\partial \zeta}(\mathcal{H}(\zeta)\varepsilon(\zeta ,t)) + P_0 \mathcal{H}(\zeta)\varepsilon(\zeta,t),
\end{split}
\label{pde-pHs}
\end{equation}
where $\varepsilon(\zeta,t) \in\mathbb{R}^n$ for $\zeta\in[a,b]$ and $t\geq 0$. In addition, $P_1 = P_1^T \in \mathbb{R}^{n\times n}$ is invertible, $P_0 = -P_0^T \in \mathbb{R}^{n\times n}$, and $\mathcal{H}\in L^{\infty}([a,b];\mathbb{R}^{n\times n})$ is symmetric and satisfies $mI\leq \mathcal{H}(\zeta)$ for all $\zeta\in [a,b]$ and some constant $m>0$. The state space $\pazocal{X}:=L^2([a,b];\mathbb{R}^n)$ is endowed with the energy inner product $\psh{\varepsilon_1}{\varepsilon_2}_{\pazocal{X}} = \psh{\varepsilon_1}{\mathcal{H}\varepsilon_2}_{L^2} = \int_a^b\varepsilon_1(\zeta)^T\mathcal{H}(\zeta)\varepsilon_2(\zeta)d\zeta$, for all $\varepsilon_1, \varepsilon_2\in\pazocal{X}$. The energy associated to (\ref{pde-pHs}) is given by $E(t) = \frac{1}{2}\Vert\varepsilon(t)\Vert_\pazocal{X}^2$.

The boundary ports denoted by $f_\partial$ and $e_\partial$ are given by 
\begin{equation}
\left(
\begin{smallmatrix}
f_\partial(t)\\
e_\partial(t)
\end{smallmatrix}\right) = \frac{1}{\sqrt{2}} \left(\begin{smallmatrix}
P_1 & -P_1\\
I & I
\end{smallmatrix}\right)
\left(\begin{smallmatrix}
(\mathcal{H} \varepsilon(t))(b)\\
(\mathcal{H} \varepsilon(t))(a)
\end{smallmatrix}\right) =: R_0\left(\begin{smallmatrix}
(\mathcal{H} \varepsilon(t))(b)\\
(\mathcal{H} \varepsilon(t))(a)
\end{smallmatrix}\right)
\end{equation}
and represent a linear combination of the restriction at the boundary variables. Note that the notation $(\mathcal{H}\varepsilon(t))(a) := \mathcal{H}(a)\varepsilon(a,t)$ has been used.
We complete the PDE (\ref{pde-pHs}) with the following homogeneous boundary conditions
\begin{align}
0 = W_B \left[
\begin{matrix}
f_\partial(t) \\
e_\partial(t)
\end{matrix}
\right], 
\label{boundary}
\end{align}
where $W_B
\in \mathbb{R}^{n\times 2n}$. It can be easily seen that the PDE (\ref{pde-pHs}) together with the boundary conditions (\ref{boundary}) admits the abstract representation $\dot{\varepsilon}(t) = \mathcal{A}\varepsilon(t), \varepsilon(0) = \varepsilon_0\in\mathcal{X}$ where the linear (unbounded) operator $\mathcal{A}$ is defined by
\begin{equation}
\mathcal{A}\varepsilon := P_1 \frac{d}{d\zeta}(\mathcal{H} \varepsilon) + P_0 \mathcal{H} \varepsilon
\label{operatorJ}
\end{equation} 
for $\varepsilon$ on the domain 
\begin{equation}
D(\mathcal{A})= \left\lbrace \varepsilon\in \pazocal{X} : \mathcal{H}\varepsilon\in H^{1}([a,b];\mathbb{R}^n), W_{B}\left[
\begin{matrix}
f_{\partial}\\
e_{\partial}
\end{matrix}
\right] =0 \right\rbrace.
\label{domainJ}
\end{equation}
Before introducing the concept of Dirac structure for (\ref{pde-pHs}) and (\ref{boundary}), let us recall the following two results from \cite{zwart}.
\begin{lemma}\label{LemmaPSH}
Consider the operator $\mathcal{A}$ defined by (\ref{operatorJ}) with domain (\ref{domainJ}). Then the following result holds:
\begin{equation}
    \psh{\mathcal{A}\varepsilon}{\varepsilon}_{\pazocal{X}} + \psh{\varepsilon}{\mathcal{A}\varepsilon}_{\pazocal{X}} = 2f_\partial^Te_\partial.
    \label{relationports}
\end{equation}
\end{lemma}

\begin{thm}\label{ThmDeterministic}
Let $W_B$ be a $n\times 2n$ real matrix. Then the operator $\mathcal{A}$ defined in (\ref{operatorJ}) on the domain (\ref{domainJ}) generates a contraction $C_0$-semigroup of bounded linear operators if and only if $W_B$ is full rank and satisfies $W_B \Sigma W_B^T \geq 0$, with $\Sigma := \left(\begin{smallmatrix}0 & I\\ I & 0\end{smallmatrix}\right)$. Furthermore, the energy balance equation
\begin{equation}
\frac{dE(t)}{dt} = f_\partial^T(t) e_\partial(t)
\label{EnergyEquation}
\end{equation}
holds.
\end{thm}
\begin{pf}
For the fact that the operator $\mathcal{A}$ generates a contraction $C_0-$semigroup, we refer to \cite[Theorem 7.2.4]{zwart}. Now observe that
\begin{align*}
    \frac{dE}{dt}(t) = \frac{1}{2}\left(\psh{\mathcal{A}\varepsilon(t)}{\varepsilon(t)}_\mathcal{X} + \psh{\varepsilon(t)}{\mathcal{A}\varepsilon(t)}_\mathcal{X}\right) = f_\partial^T(t)e_\partial(t)
\end{align*}
according to Lemma \ref{LemmaPSH}, which ends the proof.$\qed$
\end{pf}
The energy balance equation (\ref{EnergyEquation}) allows us to introduce the notion of Dirac structure for deterministic PHSs. In that way let us define the flow and the effort spaces as follows 
\begin{equation}
\mathcal{F}=\mathcal{E}=\pazocal{X}\times \mathbb{R}^n.
\end{equation}
The pairing $\left\langle \cdot,\cdot\right\rangle_{+}$ is then given by 
\begin{align}
&\psh{(f^1, f^1_\partial, e^1, e_\partial^1)}{(f^2, f^2_\partial, e^2, e_\partial^2)}_+\nonumber\\
&= \psh{e^1}{f^2}_{L^2}+ \psh{e^2}{f^1}_{L^2} - \psh{e^1_\partial}{f^2_\partial} - \psh{e_\partial^2}{f_\partial^1}
\label{innerproduct}
\end{align}
for any $(f^i, f^i_\partial, e^i, e^i_\partial)\in \mathcal{F}\times\mathcal{E}, i=1,2$. According to (\ref{SymmetricPairing}) the linear mapping $j:\mathcal{F}\to\mathcal{E}$ is given by\footnote{It is easy to see that $j = j^{-1}$ in that case. Moreover, we consider that the space $\mathcal{X}$ is endowed with the classical inner product $\langle\cdot,\cdot\rangle_{L^2}$ (which is equivalent to $\langle\cdot,\cdot\rangle_\mathcal{X}$) when defining the pairing $\langle\cdot,\cdot\rangle_+$.} $j = \left(\begin{smallmatrix}I_\mathcal{X} & 0\\ 0 & -I_{\mathbb{R}^n}\end{smallmatrix}\right)$. By considering the flow and the effort variables as \begin{equation}
f_\varepsilon = \frac{\partial \varepsilon}{\partial t} \qquad \text{and} \qquad e_\varepsilon=\mathcal{H}\varepsilon,
\label{PortsPHS}
\end{equation}
a linear first order port-Hamiltonian system is then described by 
\begin{equation*}
\left\lbrace \varepsilon(\cdot,t) | \left(\begin{matrix}
f_\varepsilon\\
f_\partial\\
e_\varepsilon\\
e_\partial
\end{matrix}\right)\in \mathcal{D} \right\rbrace,
\end{equation*}
where $\mathcal{D}$ is given as
\begin{align}
\mathcal{D}&= \left\lbrace \left(\begin{matrix}
f_\varepsilon\\
f_\partial\\
e_\varepsilon\\
e_\partial
\end{matrix}\right) \in \mathcal{F}\times \mathcal{E} | e_\varepsilon \in H^1([a,b];\mathbb{R}^n), f_\varepsilon=\mathcal{J}e_\varepsilon,\right.\nonumber\\
&\hspace{2cm}\left.\left(\begin{matrix}
f_\partial\\
e_\partial 
\end{matrix}\right) = R_0 \left(\begin{matrix}
(e_\varepsilon)(b)\\
(e_\varepsilon)(a)
\end{matrix}\right)
\right\rbrace,
\label{DiracDeterministic}
\end{align}
in which the operator $\mathcal{J}$ is defined via the operator $\mathcal{A}$ as $\mathcal{J}(\mathcal{H}\varepsilon) := \mathcal{A}\varepsilon$ for all $\varepsilon\in D(\mathcal{A})$.

\section{Stochastic port-Hamiltonian systems}\label{sec:SPHSs}

This section is devoted to the extension of the notion of Dirac structure to SPHSs. In particular, new ports are considered in order to take the stochastic effects into account, and a specific structure built from (\ref{DiracDeterministic}) is shown to be a Dirac structure for a class of SPHSs.

Let $(\Omega, \mathcal{F}, \mathbb{F}, \mathbb{P})$ be a complete filtered probability space, where $\Omega$
denotes the sample space, $\mathcal{F}$ denotes a $\sigma$-algebra, $\mathbb{F}:=(\mathcal{F}_t)_{t\geq 0}$ is a normal filtration, and $\mathbb{P}$ defines a probability measure. The class of stochastic port-Hamiltonian systems is governed by the following stochastic partial differential equation (SPDE)
\begin{align}
\dfrac{\partial \varepsilon}{\partial t}(\zeta ,t) = P_1 &\dfrac{\partial}{\partial \zeta}(\mathcal{H}(\zeta)\varepsilon(\zeta ,t)) + P_0 \mathcal{H}(\zeta)\varepsilon(\zeta,t)\nonumber\\
&+ H(\varepsilon(\zeta,t))(\dot{w}(t))(\zeta)
\label{spde}
\end{align}
with the boundary conditions (\ref{boundary}).  Let $Z$ be a Hilbert space. The noise process $(\dot{w}(t))_{t\geq 0}$ is a $Z-$valued Gaussian white noise process with intensity operator $H\in\pazocal{L}(\pazocal{X},\pazocal{L}(Z,\pazocal{X}))$. The SPDE (\ref{spde}) can be rewritten as a SDE on $\mathcal{X}$ of the form 
\begin{equation}
\delta\varepsilon(t) = \mathcal{A} \varepsilon(t) \delta t + H(\varepsilon(t)) \delta w(t), \varepsilon(0) = \varepsilon_0,
\label{SDE}
\end{equation}
where the operator $\mathcal{A}$ is given in (\ref{operatorJ}) with domain (\ref{domainJ}). Note that $(w(t))$ is a $Z-$valued Wiener process with covariance operator $Q\in\pazocal{L}(Z)$. We further assume that $Q$ is nonnegative, self-adjoint and of trace class, i.e. $\text{Tr}[Q] < \infty$. The differential term $H(\varepsilon(t))\delta w(t)$ has to be understood under the Stratonovich definition of a stochastic integral. For further details on the Stratonovich stochastic integral on Hilbert spaces, we refer to \cite{duan2014} among others. It is worth pointing out that the Stratonovich integral satisfies the standard rules of the chain calculus.
The well-posedness of the SDE (\ref{SDE}) in terms of existence and uniqueness of a mild solution is extensively studied in \cite{LamolineIEEE} and \cite{LamolineEjc}. By a mild solution of (\ref{SDE}), we mean a $\mathbb{F}$-adapted and mean-square continuous solution of the integral form of (\ref{SDE}), i.e. a solution $\varepsilon(t)$ which satisfies
\begin{align*}
    \varepsilon(t) = T(t)\varepsilon_0 + \int_0^t T(t-s)H(\varepsilon(s))\delta w(s),
\end{align*}
where $(T(t))_{t\geq 0}$ is the $C_0-$semigroup\footnote{We assume being in the conditions of Theorem \ref{ThmDeterministic}.} whose operator $\mathcal{A}$ is the infinitesimal generator and $\varepsilon_0\in\pazocal{X}$ denotes the initial condition. The power-balance equation associated to (\ref{spde}) can be expressed as follows
\begin{equation*}
\begin{split}
\delta E(\varepsilon(t)) = \delta f_\partial^T(t) &e_\partial(t)+ \psh{H^*(\varepsilon(t)) \mathcal{H}\varepsilon(t)}{\delta w(t)}_Z,
\end{split}
\end{equation*} 
which is equivalent to
\begin{equation}
0 = -\psh{\delta f_\varepsilon(t)}{e_\varepsilon(t)}_{L^2} + \delta f_\partial^T(t) e_\partial(t) + \psh{\delta w(t)}{e_w(t)}_Z,
\end{equation}
where the ports $\delta f_\varepsilon =: \delta\varepsilon$ and $e_\varepsilon = \mathcal{H}\varepsilon$ while the new ports (noise ports) due to the stochastic nature of (\ref{spde}) are defined as $\delta f_w:=\delta w(t)$ and $e_w(t) = H^*e_\varepsilon$. Note that the boundary ports $\delta f_\partial$ and $e_\partial$ are now expressed as $\delta f_\partial = \frac{1}{\sqrt{2}} (P_1 (e_\varepsilon)(b) - P_1 (e_\varepsilon)(a))\delta t$ and $e_\partial = \frac{1}{\sqrt{2}}(e_\varepsilon(b) - e_\varepsilon(a))$, respectively. In that way the pairing (\ref{innerproduct}) is extended as follows
\begin{align}
&\psh{(\delta f^1_\varepsilon, \delta f^1_\partial,\delta f_w^1, e^1_\varepsilon, e_\partial^1,e_w^1)}{(\delta f^2_\varepsilon, \delta f^2_\partial, \delta f_w^2, e^2_\varepsilon, e_\partial^2, e_w^2)}_+\nonumber\\
&= \psh{e^1_\varepsilon}{\delta f^2_\varepsilon}_{L^2} + \psh{e^2_\varepsilon}{\delta f^1_\varepsilon}_{L^2} - \psh{e^1_\partial}{\delta f^2_\partial}_{\mathbb{R}^n}\nonumber\\
&- \psh{e_\partial^2}{\delta f_\partial^1}_{\mathbb{R}^n} - \psh{e_w^1}{\delta f_w^2}_Z - \psh{e_w^2}{\delta f_w^1}_Z.\label{innerproductsto}
\end{align}
Observe that $e_w(t)$ represents the power-conjugated effort coupled to the stochastic input $w(t)$. As it will be shown in Theorem \ref{MainTheorem}, it enables to preserve Dirac structure of port-Hamiltonian systems when subject to stochastic disturbances, see \citep{LamolineEjc}.
\begin{rem}
We stress that the notations $\delta f_\varepsilon, \delta f_\partial$ and $\delta f_w$ are used for the flows as they are defined from infinitesimal variations resulting from $\delta t$ and $\delta w(t)$. 
\end{rem}
In order to write SPHSs as Dirac structures, we complete the flow and the efforts spaces in the following way
\begin{equation}
\mathcal{F} = \mathcal{E} := \pazocal{X}\times Z\times \mathbb{R}^n.
\label{FlowEffortSpacesStoch}
\end{equation}
Once more, the comparison can be made with (\ref{SymmetricPairing}) where the invertible linear map $j:\mathcal{F}\to\mathcal{E}$ is expressed as $j = \left(\begin{smallmatrix}I_\mathcal{X} & 0 & 0\\ 0 & -I_Z & 0\\ 0 & 0 & -I_{\mathbb{R}^n}\end{smallmatrix}\right)$. As a result of the pairing (\ref{innerproductsto}), let us consider the following structure for SPHSs described by (\ref{spde})
\begin{align}
&\mathcal{D}= \left\lbrace\left(\begin{matrix}
\delta f_\varepsilon\\
\delta f_w\\
\delta f_\partial\\
e_\varepsilon\\
e_w\\
e_\partial
\end{matrix}\right) \in \mathcal{F}\times \mathcal{E} | e_\varepsilon \in H^1([a,b];\mathbb{R}^n),\right.\nonumber\\
&\delta f_\varepsilon=\mathcal{J}e_\varepsilon \delta t + H \delta f_w, \delta f_\partial  = \frac{1}{\sqrt{2}} (P_1 (e_\varepsilon)(b) - P_1 (e_\varepsilon)(a))\delta t,\nonumber\\
&\left.e_\partial = \frac{1}{\sqrt{2}}(e_\varepsilon(b) - e_\varepsilon(a)) , e_w = H^* e_\varepsilon\right\rbrace.
\label{diracstochst}
\end{align}

Note that the boundary flow and effort variables, $\delta f_\partial$ and $e_\partial$ respectively, can be written in a more compact form as
\begin{equation}
\left(\begin{matrix}
\delta f_\partial\\
e_\partial \delta t
\end{matrix}\right) = R_0\left(\begin{matrix}(e_\varepsilon)(b)\\
(e_\varepsilon)(a)\end{matrix}\right)\delta t.
\label{BoundPortsCompact}
\end{equation}
We can now prove the main result, which states that $\mathcal{D}$ given by (\ref{diracstochst}) forms a Dirac structure as defined in Definition \ref{diracdef}.
\begin{thm}\label{MainTheorem}
The subspace $\mathcal{D}$ of $\mathcal{B}$ given by (\ref{diracstochst}) is a Dirac structure. \end{thm}
\begin{pf}
We first prove that $\mathcal{D}\subset \mathcal{D}^\perp$. This is equivalent to the canonical product $\psh{b}{b}_+$ being set to zero for any $b\in\mathcal{D}$. Let us consider $(\delta f_\varepsilon,\delta f_w,\delta f_\partial,e_\varepsilon,e_w,e_\partial)\in\mathcal{F}\times\mathcal{E}$. From (\ref{innerproductsto}), we get that 
\begin{align*}
&{\psh{(\delta f_\varepsilon,\delta f_w,\delta f_\partial, e_\varepsilon,e_w,e_\partial)}{(\delta f_\varepsilon,\delta f_w,\delta f_\partial, e_\varepsilon,e_w,e_\partial)}_+}\\
&= \psh{\delta f_\varepsilon}{e_\varepsilon}_{L^2} - \psh{\delta f_w}{e_w}_Z -\psh{\delta f_\partial}{e_\partial}_{\mathbb{R}^n}  \\&+ \psh{\delta f_\varepsilon}{e_\varepsilon}_{L^2}  - \psh{\delta f_w}{e_w}_Z  - \psh{\delta f_\partial}{e_\partial}_{\mathbb{R}^n}.
\end{align*}
Moreover, since $\delta f_\varepsilon(t)=\mathcal{J}e_\varepsilon(t)\delta t + H(\varepsilon(t)) \delta f_w(t)$, we obtain that
\begin{align*}
&\psh{(\delta f_\varepsilon,\delta f_w,\delta f_\partial, e_\varepsilon,e_w,e_\partial)}{(\delta f_\varepsilon,\delta f_w,\delta f_\partial, e_\varepsilon,e_w,e_\partial)}_+\\
&= \psh{\mathcal{J}e_\varepsilon(t) \delta t + H(\varepsilon(t)) \delta f_w(t)}{e_\varepsilon(t)}_{L^2} - \psh{\delta f_\partial(t)}{e_\partial(t)}_{\mathbb{R}^n}\\
&- \psh{\delta f_w(t)}{H^*(\varepsilon(t))e_\varepsilon(t)}_Z\\
&+ \psh{\mathcal{J}e_\varepsilon(t) \delta t + H(\varepsilon(t)) \delta f_w(t)}{e_\varepsilon(t)}_{L^2}\\
&- \psh{\delta f_\partial(t)}{e_\partial(t)}_{\mathbb{R}^n} - \psh{\delta f_w(t)}{H^*(\varepsilon(t))e_\varepsilon(t)}_Z.
\end{align*}
By using (\ref{diracstochst}), there holds
\begin{align*}
&\psh{(\delta f_\varepsilon,\delta f_w,\delta f_\partial, e_\varepsilon,e_w,e_\partial)}{(\delta f_\varepsilon,\delta f_w,\delta f_\partial, e_\varepsilon,e_w,e_\partial)}_+\\
&= \psh{\mathcal{J}e_\varepsilon(t)}{e_\varepsilon(t)}_{L^2} \delta t  + \psh{\mathcal{J}e_\varepsilon(t)}{e_\varepsilon(t)}_{L^2} \delta t\\
&- 2 \psh{\delta f_\partial(t)}{e_\partial(t)}_{\mathbb{R}^n} =0.
\end{align*}
In order to prove that $\mathcal{D}^\perp\subset \mathcal{D}$, let us pick any $(\delta f_\varepsilon, \delta f_w, \delta f_\partial, e_\varepsilon, e_w, e_\partial) \in \mathcal{D}^\perp$ and $(\delta \tilde{f}_\varepsilon, \delta \tilde{f}_w, \delta \tilde{f}_\partial, \tilde{e}_\varepsilon, \tilde{e}_w, \tilde{e}_\partial) \in \mathcal{D}$. 
By orthogonality, we have that\footnote{The time dependency of the variables has been willingly omitted for the sake of readability.}
\begin{align}
0 &= \psh{(\delta f_\varepsilon, \delta f_w,\delta f_\partial, e_\varepsilon, e_w, e_\partial)}{(\delta \tilde{f}_\varepsilon, \delta \tilde{f}_w,\delta \tilde{f}_\partial, \tilde{e}_\varepsilon, \tilde{e}_w,\tilde{e}_\partial)}_+\nonumber\\
&=- \psh{\delta f_\varepsilon}{\tilde{e}_\varepsilon}_{L^2} +  \psh{\delta f_\partial}{\tilde{e}_\partial}_{\mathbb{R}^n}  + \psh{\delta f_w}{\tilde{e}_w}_Z\nonumber\\
&- \psh{\delta \tilde{f}_\varepsilon}{e_\varepsilon}_{L^2} + \psh{\delta \tilde{f}_\partial}{e_\partial}_{\mathbb{R}^n} +  \psh{\delta \tilde{f}_w}{e_w}_Z.\label{maineqpr}
\end{align}
Step 1. By setting $\tilde{e}_\varepsilon = 0$, it implies that $\delta \tilde{f}_\partial = \tilde{e}_\partial = \tilde{e}_w = 0$. Hence, 
\begin{align*}
0 &= - \psh{
H(\varepsilon(t))\delta \tilde{f}_w}{e_\varepsilon}_{L^2} + \psh{\delta \tilde{f}_w}{e_w}_Z\\
&=  \psh{\delta \tilde{f}_w}{-H^*(\varepsilon(t))e_\varepsilon+e_w}_Z.
\end{align*} 
We can deduce that $e_w = H^*(\varepsilon(t)) e_\varepsilon$. \\
Step 2. Let us now choose $\tilde{e}_\varepsilon\in H^1([a,b];\mathbb{R}^n)$ with compact support strictly included in $(a,b)$, which entails that $\tilde{e}_\varepsilon$ is zero in $a$ and $b$. We also set $\delta \tilde{f}_w = 0$. It is easy to see that $(\mathcal{J}\tilde{e}_\varepsilon \delta t, 0, 0, \tilde{e}_\varepsilon, \tilde{e}_w, 0) \in \mathcal{D}$. Therefore, (\ref{maineqpr}) becomes 
\begin{align*}
0 &=  {\psh{(\delta f_\varepsilon, \delta f_w,\delta f_\partial, e_\varepsilon, e_w, e_\partial)}{(\mathcal{J}\tilde{e}_\varepsilon \delta t, 0, 0, \tilde{e}_\varepsilon, \tilde{e}_w, 0)}_+}\\
&= -\psh{\delta f_\varepsilon}{\tilde{e}_\varepsilon}_{L^2} + \psh{H(\varepsilon) \delta f_w}{\tilde{e}_\varepsilon}_{L^2} - \psh{e_\varepsilon}{\mathcal{J}\tilde{e}_\varepsilon}_{L^2}\delta t.
\end{align*} 
An integration by parts on the term $\psh{e_\varepsilon}{\mathcal{J}\tilde{e}_\varepsilon}_{L^2}:=\int_a^b e_\varepsilon(\zeta)^T \mathcal{J}\tilde{e}_\varepsilon(\zeta) d\zeta$ yields that 
\begin{equation*}
\int_a^b e_\varepsilon(\zeta)^T \mathcal{J}\tilde{e}_\varepsilon(\zeta) d\zeta = [ e_\varepsilon^T P_1 \tilde{e}_\varepsilon]_a^b - \psh{\mathcal{J}e_\varepsilon}{\tilde{e}_\varepsilon}_{L^2}.
\end{equation*}
Therefore, we obtain $\delta f_\varepsilon(t)= \mathcal{J}e_\varepsilon(t)\delta t + H(\varepsilon(t)) \delta f_w(t)$.\\
Step 3. Let us take now $\delta \tilde{f}_w=\tilde{e}_w=0$ such that $(\mathcal{J}\tilde{e}_\varepsilon \delta t, 0, \delta \tilde{f}_\partial, \tilde{e}_\varepsilon, 0, \tilde{e}_\partial) \in \mathcal{D}$. From (\ref{maineqpr}), we have
\begin{align*}
&0 = \psh{\delta f_\varepsilon}{\tilde{e}_\varepsilon}_{L^2} + \psh{\delta \tilde{f}_\varepsilon}{e_\varepsilon}_{L^2} - \psh{\delta f_\partial}{\tilde{e}_\partial}_{\mathbb{R}^n} - \psh{\delta \tilde{f}_\partial}{e_\partial}_{\mathbb{R}^n}\\
&= \psh{\mathcal{J}e_\varepsilon}{\tilde{e}_\varepsilon}_{L^2}\delta t + \psh{\mathcal{J}\tilde{e}_\varepsilon}{e_\varepsilon}_{L^2}\delta t\\
&\hspace{4.3cm}- \psh{\delta f_\partial}{\tilde{e}_\partial}_{\mathbb{R}^n} - \psh{\delta \tilde{f}_\partial}{e_\partial}_{\mathbb{R}^n}\\
&=[ e_\varepsilon^T P_1 \tilde{e}_\varepsilon]_a^b \delta t - \psh{\delta f_\partial}{\tilde{e}_\partial}_{\mathbb{R}^n} - \psh{\delta \tilde{f}_\partial}{e_\partial}_{\mathbb{R}^n}\\
&=\left[
\begin{matrix}
\delta \tilde{f}_\partial \\
\tilde{e}_\partial \delta t  
\end{matrix} 
\right]^T \Sigma R_0 \left[ 
\begin{matrix}
e_\varepsilon(b)\\
e_\varepsilon(a)
\end{matrix}
\right] - \psh{\delta f_\partial}{\tilde{e}_\partial}_{\mathbb{R}^n} - \psh{\delta \tilde{f}_\partial}{e_\partial}_{\mathbb{R}^n}\\
&=\left[
\begin{matrix}
\tilde{e}_\partial \delta t\\
\delta \tilde{f}_\partial  
\end{matrix} 
\right]^T R_0 
\left[ 
\begin{matrix}
(e_\varepsilon)(b)\\
(e_\varepsilon)(a)
\end{matrix}
\right] - \psh{\delta f_\partial}{\tilde{e}_\partial}_{\mathbb{R}^n} - \psh{\delta \tilde{f}_\partial}{e_\partial}_{\mathbb{R}^n}\\
&= \left[
\begin{matrix}
\tilde{e}_\partial \delta t\\
\delta \tilde{f}_\partial  
\end{matrix} 
\right]^T
\left( 
R_0 
\left[ 
\begin{matrix}
(e_\varepsilon)(b)\\
(e_\varepsilon)(a)
\end{matrix}
\right] - \left[
\begin{matrix}
\delta f_\partial\frac{1}{\delta t} \\
e_\partial  
\end{matrix} 
\right]
\right).
\end{align*}
Since the above equality has to hold for all $\tilde{e}_\partial$ and $\delta \tilde{f}_\partial$, we deduce that 
\begin{equation*}
\left(
\begin{matrix}
\delta f_\partial\\
e_\partial \delta t
\end{matrix}
\right) = R_0 \left(\begin{matrix} 
(\mathcal{H}\varepsilon(t))(b)\\
(\mathcal{H}\varepsilon(t))(a)
\end{matrix}
\right)\delta t.
\end{equation*}
This proves that  $\mathcal{D}^\perp\subset \mathcal{D}$, which completes the proof.$\qed$ 
\end{pf}

\begin{rem}
Dissipative elements have not been considered to focus on proof arguments regarding the noise elements. Note that dissipative elements could be added independently to the Dirac structure. Theorem \ref{MainTheorem} would then readily extend, see \cite[Theorem 6.5]{Villegas}. 
\end{rem}

\section{Illustration: boundary control as interconnection of stochastic Dirac structures}\label{InterSection}
In this section, we illustrate a central feature of Dirac structures, their ability to be interconnected between each other, under appropriate assumptions. We investigate the interconnection of a stochastic Dirac structure of the form (\ref{diracstochst}) with an arbitrary Dirac structure that shares the same space of interconnection. Therefore, the notions of split Tellegen and split Dirac structures are introduced, see \cite[Definitions 3.1]{Kurula}.
\begin{defn}
Let us suppose that the flow and the effort spaces may be decomposed as $\mathcal{F} = \mathcal{F}_1\times \mathcal{F}_2$ and $\mathcal{E} = \mathcal{E}_1\times\mathcal{E}_2$, respectively. Furthermore, assume that $j_i:\mathcal{F}_i\to\mathcal{E}_i, i=1, 2$ are unitary linear and invertible mappings. A linear subspace $\mathcal{D}\subset\mathcal{B}$ is called a split Tellegen (split Dirac) structure if it is a Tellegen (Dirac) structure in the sense of Definition \ref{diracdef} with $j = \left(\begin{smallmatrix}j_1 & 0\\ 0 & j_2\end{smallmatrix}\right)$.
\end{defn}
We are now in position to explicit what we mean by the interconnection of two (stochastic) Dirac structures. Therefore, let us consider $\mathcal{D}^A\subset \mathcal{B} = (\mathcal{F}_1\times\mathcal{F}_2)\times (\mathcal{E}_1\times\mathcal{E}_2)$ being a split Dirac structure of the form (\ref{diracstochst}) and $\mathcal{D}^B\subset \mathcal{B} = (\mathcal{F}_3\times\mathcal{F}_2)\times (\mathcal{E}_3\times\mathcal{E}_2)$ being an arbitrary split Dirac structure where $\mathcal{F}_3$ and $\mathcal{E}_3$ are Hilbert spaces. As $\mathcal{D}^A$ is of the form (\ref{diracstochst}), we define the spaces $\mathcal{F}_1, \mathcal{E}_1$ and $\mathcal{F}_2, \mathcal{E}_2$ as $\mathcal{F}_1 = \mathcal{X}\times Z\times\mathbb{R}^{n-p} = \mathcal{E}_1$ and $\mathcal{F}_2 = \mathbb{R}^p = \mathcal{E}_2$, $1\leq p\leq n$, respectively. With the split Dirac structures $\mathcal{D}^A$ and $\mathcal{D}^B$, we define the unitary operators $j_1: \mathcal{F}_1\to\mathcal{E}_1, j_2:\mathcal{F}_2\to\mathcal{E}_2$ and $j_3:\mathcal{F}_3\to\mathcal{E}_2$ where $j_1$ and $j_2$ are given by\footnote{From this definition of $j_1$ and $j_2$, there holds $j_1=j_1^{-1}$ and $j_2 = j_2^{-1}$.}
\begin{align*}
    j_1 = \left(\begin{smallmatrix}
        I_\mathcal{X} & 0 & 0\\
        0 & -I_Z & 0\\
        0 & 0 & -I_{\mathbb{R}^{n-p}}
    \end{smallmatrix}\right),\,\, j_2 = -I_{\mathbb{R}^p}.
\end{align*}
The proposed interconnection between $\mathcal{D}^A$ and $\mathcal{D}^B$ is inspired by \cite[Definition 3.2]{Kurula} in which first-order PHSs are considered. It is performed via some of the boundary ports through the space $\mathcal{F}_2\times \mathcal{E}_2$. The dimension of $\mathcal{F}_2$ gives the number of boundary ports used for the interconnection.
\begin{figure*}
    \centering
    \includegraphics[scale=1,trim=0cm 1cm 0cm 1.2cm]{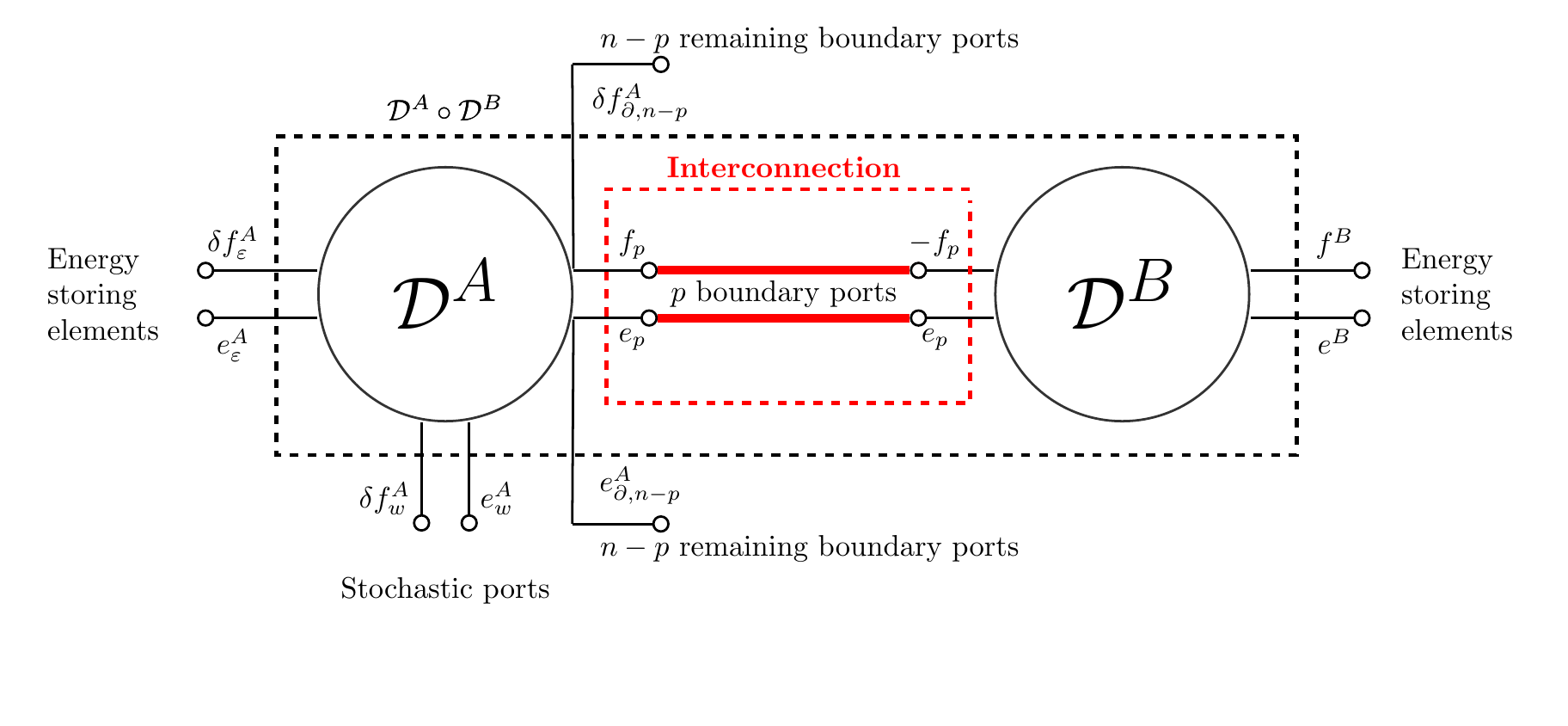}
    \caption{Interconnection of the Dirac structures $\mathcal{D}^A$ and $\mathcal{D}^B$. The resulting structure is denoted by $\mathcal{D}^A\circ\mathcal{D}^B$. The interconnection is performed at the boundary thanks to the boundary ports $f_p$ and $e_p$.}
    \label{fig:Interconnection}
\end{figure*}
\begin{defn}\label{InterconnectionDefinition}
The composition of $\mathcal{D}^A$ and $\mathcal{D}^B$ is denoted $\mathcal{D}^A\circ \mathcal{D}^B\subset (\mathcal{F}_1\times\mathcal{F}_3)\times (\mathcal{E}_1\times\mathcal{E}_3)$ and is defined as
\begin{align}
\mathcal{D}^A\circ\mathcal{D}^B := \left\{\left(\begin{smallmatrix}\delta f^A\\ f^B\\ e^A\\ e^B\end{smallmatrix}\right)\in (\mathcal{F}_1\times\mathcal{F}_3)\times (\mathcal{E}_1\times\mathcal{E}_3)\right.\nonumber\\
\left. \exists f_p\in\mathcal{F}_2, e_p\in\mathcal{E}_2, \left(\begin{smallmatrix}\delta f^A\\ f_p\\ e^A\\ e_p\end{smallmatrix}\right)\in\mathcal{D}^A, \left(\begin{smallmatrix}f^B\\ -f_p\\ e^B\\ e_p\end{smallmatrix}\right)\in\mathcal{D}^B\right\}\label{InterconnectedStructure}
\end{align}
with\footnote{This definition of $j$ entails that $j^{-1} = \left(\begin{smallmatrix}j_1 & 0\\ 0 & j_3^{-1}\end{smallmatrix}\right)$.} $j = \left(\begin{smallmatrix}j_1 & 0\\0 & j_3\end{smallmatrix}\right)$, where\footnote{The notation $\delta f_{\partial,n-p}$ is used to emphasize the fact that the vector $\delta f_{\partial,n-p}$ is of size $n-p$.} $\delta f^A = \left(\begin{smallmatrix}\delta f_\varepsilon^A, & \delta f_w^A, & \delta f^A_{\partial,n-p}\end{smallmatrix}\right), e^A = \left(\begin{smallmatrix} e_\varepsilon^A, & e_w^A, & e^A_{\partial,n-p}\end{smallmatrix}\right)$.
\end{defn}
An illustration of the proposed interconnection is given in Figure \ref{fig:Interconnection}. In that way, it is easy to see that the pairing whose $\mathcal{D}^A\circ \mathcal{D}^B$ is equipped with, denoted $\langle\cdot,\cdot\rangle_\dagger$, is expressed as
\begin{align}
&\left\langle\left(\begin{smallmatrix}\delta f^A\\ f^B\\ e^A\\ e^B\end{smallmatrix}\right),\left(\begin{smallmatrix}\delta \tilde{f}^A\\ \tilde{f}^B\\ \tilde{e}^A\\ \tilde{e}^B\end{smallmatrix}\right)\right\rangle_\dagger = \left\langle\left(\begin{smallmatrix}\delta f_\varepsilon^A\\ \delta f_w^A\\ \delta f^A_{\partial,n-p}\\ f^B\end{smallmatrix}\right), j^{-1}\left(\begin{smallmatrix}\tilde{e}_\varepsilon^A\\ \tilde{e}_w^A\\ \tilde{e}^A_{\partial,n-p} \\\tilde{e}^B\end{smallmatrix}\right)\right\rangle\nonumber\\
&+ \left\langle\left(\begin{smallmatrix}e_\varepsilon^A\\ e_w^A\\ e^A_{\partial,n-p}\\ e^B\end{smallmatrix}\right), j\left(\begin{smallmatrix}\delta \tilde{f}_\varepsilon^A\\ \delta \tilde{f}_w^A\\ \delta \tilde{f}^A_{\partial,n-p}\\ \tilde{f}^B\end{smallmatrix}\right)\right\rangle = \langle \delta f_\varepsilon^A,\tilde{e}_\varepsilon^A\rangle - \langle \delta f_w^A,\tilde{e}_w^A\rangle\nonumber\\
& - \langle \delta f^A_{\partial,n-p}, \tilde{e}^A_{\partial,n-p}\rangle + \langle f^B, j_3^{-1}\tilde{e}^B\rangle + \langle e_\varepsilon^A, \delta \tilde{f}_\varepsilon^A\rangle\nonumber\\
&- \langle e_w^A,\delta\tilde{f}_w^A\rangle - \langle e^A_{\partial,n-p}, \delta \tilde{f}^A_{\partial,n-p}\rangle+\langle e^B, j_3\tilde{f}^B\rangle.\label{ExtendedPairingInter}
\end{align}
We shall now focus on the nature of the structure introduced in (\ref{InterconnectedStructure}). First let us consider the following lemma.
\begin{lemma}\label{LemmaTellegen}
The structure $\mathcal{D}^A\circ\mathcal{D}^B$ defined in (\ref{InterconnectedStructure}) has zero power, i.e. $\langle \mathfrak{d},\mathfrak{d}\rangle_\dagger = 0$ for any $\mathfrak{d}\in\mathcal{D}^A\circ\mathcal{D}^B$.
\end{lemma}
\begin{pf}
Let us consider $\mathfrak{d} = \left(\begin{smallmatrix}f\\ e\end{smallmatrix}\right)\in\mathcal{D}^A\circ\mathcal{D}^B$, where $f = \left(\begin{smallmatrix}\delta f_\varepsilon^A, & \delta f_w^A, & \delta f^A_{\partial,n-p}, & f^B\end{smallmatrix}\right), e = \left(\begin{smallmatrix}e_\varepsilon^A, & e_w^A, & e^A_{\partial,n-p}, & e^B\end{smallmatrix}\right)$. Showing that $\langle \mathfrak{d},\mathfrak{d}\rangle_{(+,\circ)} = 0$ is equivalent in showing that
\begin{align}
&\langle \delta f_\varepsilon^A,e_\varepsilon^A\rangle - \langle \delta f_w^A,e_w^A\rangle - \langle \delta f^A_{\partial,n-p}, e^A_{\partial,n-p}\rangle + \langle f^B, j_3^{-1}e^B\rangle\nonumber\\
&+ \langle e_\varepsilon^A, \delta f_\varepsilon^A\rangle - \langle e_w^A,\delta f_w^A\rangle - \langle e^A_{\partial,n-p},\delta f^A_{\partial,n-p}\rangle + \langle e^B,j_3 f^B\rangle\nonumber\\
&= 0,\label{ZeroPower}
\end{align}
see (\ref{ExtendedPairingInter}). According to the definition of $\mathcal{D}^A\circ\mathcal{D}^B$, there exist $f_p, e_p \in\mathbb{R}^p$ such that
$\left(\begin{smallmatrix}\delta f_\varepsilon^A, & \delta f_w^A, & \delta f_{\partial,n-p}, & f_p, & e_\varepsilon^A, & e_w^A, & e_{\partial,n-p}, & e_p\end{smallmatrix}\right)\in\mathcal{D}^A$ and
     $\left(\begin{smallmatrix} f^B, & - f_p, & e^B, & e_p\end{smallmatrix}\right)\in\mathcal{D}^B$. As $\mathcal{D}^A$ and $\mathcal{D}^B$ are split Dirac structures, they have zero power. In particular, there holds
\begin{align}
&0=\langle \delta f_\varepsilon^A,e_\varepsilon^A\rangle - \langle \delta f_w^A,e_w^A\rangle - \langle\delta f^A_{\partial,n-p},e^A_{\partial,n-p}\rangle - \langle f_p,e_p\rangle\nonumber\\
&+ \langle e_\varepsilon^A,\delta f_\varepsilon^A\rangle - \langle e_w^A,\delta f_w^A\rangle - \langle e^A_{\partial,n-p},\delta f^A_{\partial,n-p}\rangle - \langle e_p,f_p\rangle
\label{ZeroPowerOne} 
\end{align}
and
\begin{align}
&0=\langle f^B,j_3^{-1}e^B\rangle + \langle f_p,e_p\rangle + \langle e^B,j_3 f^B\rangle + \langle e_p,f_p\rangle.
\label{ZeroPowerTwo}
\end{align}
Making the sum of (\ref{ZeroPowerOne}) and (\ref{ZeroPowerTwo}) implies (\ref{ZeroPower}).\qed 
\end{pf}
This results is also known as the fact that the interconnection of split Tellegen structures remains a Tellegen structure, when the interconnection is expressed like it is in Definition \ref{InterconnectionDefinition}, see \cite[Proposition 3.3]{Kurula}. This means that $\mathcal{D}^A\circ\mathcal{D}^B\subset (\mathcal{D}^A\circ\mathcal{D}^B)^\perp$ where the orthogonal complement "$\perp$" has to be understood with respect to the new pairing $\langle\cdot,\cdot\rangle_\dagger$. However, showing the inclusion $(\mathcal{D}^A\circ\mathcal{D}^B)^\perp\subset \mathcal{D}^A\circ\mathcal{D}^B$ poses delicate problems and comes with conditions since it depends on the structures $\mathcal{D}^A, \mathcal{D}^B$ and on the nature of the interconnection. In our setting, the following proposition holds.
\begin{prop}\label{PropInter}
The structure $\mathcal{D}^A\circ\mathcal{D}^B$ defined in (\ref{InterconnectionDefinition}) is a split Dirac structure.
\end{prop}
\begin{pf}
The inclusion $\mathcal{D}^A\circ\mathcal{D}^B\subset (\mathcal{D}^A\circ\mathcal{D}^B)^\perp$ holds from Lemma \ref{LemmaTellegen}. The other inclusion follows by the fact that the space $\mathcal{F}_2\times\mathcal{E}_2 = \mathbb{R}^p\times\mathbb{R}^p$ through which the interconnection takes place is finite-dimensional, see \cite[Corollary 3.9]{Kurula}.
\qed
\end{pf}
\begin{rem}
\begin{enumerate}
    \item As it is highlighted in \cite[Corollary 3.9]{Kurula}, the dimensionality of the interconnection plays an important role in determining whether $\mathcal{D}^A\circ\mathcal{D}^B$ is a split Dirac structure or not. No conclusion could have been possible without a finite-dimensional space of interconnection $\mathcal{F}_2\times\mathcal{E}_2$. In a more general case, one should consider the scattering operators describing each of the split Dirac structures $\mathcal{D}^A$ and $\mathcal{D}^B$, see \cite[Corollary 2.8]{Kurula}. Then, conditions on these scattering operators are proposed in \cite[Theorem 3.8]{Kurula} to ensure that $\mathcal{D}^A\circ\mathcal{D}^B$ is a split Dirac structure. This result should be of interest in the case where the interconnection is performed via Hilbert spaces-valued ports.
    \item The definition of the structure $\mathcal{D}^B$ does not exclude stochastic ports. This could be envisaged as well through the spaces $\mathcal{F}_3$ and $\mathcal{E}_3$.
    \item Without any loss of generality the problem of interconnections of multiple Dirac structures can be reduced to the problem of interconnection of two Dirac structures. 
    \end{enumerate}
\end{rem}

In terms of control practice, one usage of the Dirac structure consists in taking advantage of their nice geometric properties to design control laws for achieving certain goals via the interconnection of subsystems. Most of the current methods developed for the stabilization of infinite-dimensional port-Hamiltonian systems deal with boundary controllers, see \cite{Hashem20}. Generalization to the stochastic setting leads to even more difficulties as the noise of the plant cannot be controlled. In \cite{Haddad} noise was assumed to be vanishing at the equilibrium, which in practice would be quite restrictive in terms of configurations. Recently, theses restrictions were lift in \cite{Cordoni2022} for the generalization of energy shaping techniques using weaker concepts of Casimir function and passivity. The development of adapted control methods for infinite-dimensional SPHSs remains an uncultivated field. The authors believe that Proposition \ref{PropInter} would open the way to the development of stabilization method for infinite-dimensional SPHSs via Casimir generation or energy shaping approaches.     

To illustrate Proposition \ref{PropInter}, we study the example of a boundary controlled stochastic vibrating string described by coupled SPDE-ODE. The string is assumed to be fixed at one extremity, free at the other, and subject to some stochastic damping. Moreover, we assume that some boundary conditions are dynamic. More particularly, those are actuated by a mass-spring system. The dynamics of such a stochastic adaptive controlled system are written as
\begin{align}
\displaystyle\rho(\zeta)\frac{\partial^2 z}{\partial t^2}(\zeta,t) &= \frac{\partial}{\partial\zeta}\left(T(\zeta)\frac{\partial z}{\partial \zeta}(\zeta,t)\right) - (R_t +\dot{w}(t))\frac{\partial z}{\partial t},\vspace{0.1cm}\label{String}\\
\displaystyle\frac{\partial z}{\partial t}(a,t) &= 0,\,\, T(b)\frac{\partial z}{\partial\zeta}(b,t) = 0,\vspace{0.1cm}\label{HomogeneousBC}\\
\displaystyle \frac{1}{\sqrt{2}}\frac{\partial z}{\partial t}(b,t) &= \frac{p(t)}{m},\,\,\frac{1}{\sqrt{2}}T(a)\frac{\partial z}{\partial \zeta}(a,t) = kq(t)\label{CtrlBC}
\end{align}
where the control variables $p$ and $q$ are updated adaptively as
\begin{equation}
\begin{array}{l}
     \dot{p}(t) = -kq(t) + u(t),\\
     \dot{q}(t) = \frac{1}{m}p(t),\\
     y(t) = \frac{1}{m}p(t),
\end{array}
\label{AdaptiveInputs}
\end{equation}
with $k$ and $m$ being positive parameters. Here, $\rho(\zeta)$ and $T(\zeta)$ are the mass density and the Young modulus of the string at position $\zeta\in [a,b]$. The variable $t\geq 0$ denotes the time and $z(\zeta,t)$ is the displacement of the string at $\zeta\in[a,b]$ and $t\geq 0$. The positive frictional damping parameter $R_t$ is perturbated by a real-valued white noise $\dot{w}(t)$ with covariance $\sigma^2$. By considering $\varepsilon_1(\zeta,t) := \rho(\zeta)\frac{\partial z}{\partial t}(\zeta,t)$ and $\varepsilon_2 := \frac{\partial z}{\partial\zeta}(\zeta,t)$ as the momentum and the strain, respectively, the SDE (\ref{String}) with the homogeneous boundary conditions (\ref{HomogeneousBC}) admits the following port-Hamiltonian formulation
\begin{equation}
\left\{
\begin{array}{l}
\frac{\partial \varepsilon}{\partial t} = P_1\frac{\partial}{\partial\zeta}(\mathcal{H}\varepsilon(t)) - G_0 \mathcal{H}\varepsilon(t) + H(\varepsilon(t))(\dot{w}(t)),\\
W_B\left[\begin{smallmatrix}f_\partial(t)\\ e_\partial(t)\end{smallmatrix}\right] = \left[\begin{smallmatrix}
    0\\
    0
\end{smallmatrix}
\right],
\end{array}\right.
\label{String_PHS}
\end{equation}
with $W_B = \frac{1}{\sqrt{2}}\left[\begin{smallmatrix}1 & 0 & 0 & 1\\ 0 & -1 & 1 & 0\end{smallmatrix}\right]$ and 
\begin{align}
f_\partial(t) &= \frac{1}{\sqrt{2}}\left[\begin{smallmatrix}T(b)\frac{\partial z}{\partial \zeta}(b,t) - T(a)\frac{\partial z}{\partial \zeta}(a,t)\\ \frac{\partial z}{\partial t}(b,t) - \frac{\partial z}{\partial t}(a,t)\end{smallmatrix}\right],\nonumber\\
e_\partial(t) &= \frac{1}{\sqrt{2}}\left[\begin{smallmatrix}\frac{\partial z}{\partial t}(b,t) + \frac{\partial z}{\partial t}(a,t)\\T(b)\frac{\partial z}{\partial \zeta}(b,t) + T(a)\frac{\partial z}{\partial \zeta}(a,t)\end{smallmatrix}\right].\label{FlowEffort}
\end{align}
The Hamiltonian operator $\mathcal{H}$, the matrix $P_1$ and the matrix $G_0$ are respectively given by 
\begin{align*}
\mathcal{H}(\zeta) = \left(\begin{smallmatrix}
\frac{1}{\rho(\zeta)} & 0\\ 0 & T(\zeta)
\end{smallmatrix}\right), P_1 = \left(\begin{smallmatrix}
0 & 1\\ 1 & 0
\end{smallmatrix}\right)
\text{ and } G_0=\left(\begin{smallmatrix}
R_t & 0\\
0 & 0
\end{smallmatrix}\right).
\end{align*}
As it is made in \cite{LamolineEjc}, we set $Z = \mathbb{R}$. In that way, $H\in\mathcal{L}(\mathcal{X},\mathcal{L}(Z,\mathcal{X}))$ with $H(\varepsilon) = \left[\begin{smallmatrix}-\varepsilon_1 & 0\end{smallmatrix}\right]^T$ for all $\varepsilon\in\mathcal{X}$. Here we consider $\mathcal{X} = L^2([a,b];\mathbb{R}^2)$ as energy space. It is then easy to see that (\ref{String_PHS}) defines a Dirac structure thanks to Theorem \ref{MainTheorem}. Let us now focus on (\ref{AdaptiveInputs}) with the boundary conditions (\ref{CtrlBC}). The system (\ref{AdaptiveInputs}) may be regarded as a mass-spring system (harmonic oscillator) with $q$ and $p$ being the deviation from the zero position and the momentum, respectively. The first equation of (\ref{AdaptiveInputs}) is due to the force and the second equation is for the velocity. By defining the associated potential and kinetic energies as
\begin{equation}
    H_p = \frac{1}{2}k q^2, H_c = \frac{1}{2m} p^2,
    \label{Energies}
\end{equation}
system (\ref{AdaptiveInputs}) admits a port-Hamiltonian formulation in terms of the following Dirac structure
\begin{align}
    \mathcal{D}_c &:= \left\{\left(\begin{smallmatrix}f_1 & u & f_2 & e_1 & y & e_2\end{smallmatrix}\right)^T\in\mathbb{R}^3\times\mathbb{R}^3, f_2 = -e_1 = -y,\right.\nonumber\\
    &\left.\hspace{4cm}f_1 = e_2 - u\right\},
\end{align}
where the variables $f_1 = -\dot{p}, f_2 = -\dot{q}$ while the variables $e_1$ and $e_2$ are the derivatives of the Hamiltonian $H_p+H_c$ with respect to $p$ and $q$, respectively. The variables $u$ and $y$ are external input and output whose objective could be the stability of the closed-loop system (\ref{String}) -- (\ref{AdaptiveInputs}) for instance. As a particularity of Dirac structure, note that the power associated to $\mathcal{D}_c$ is zero, i.e. $f_1e_1 + f_2e_2 + uy = 0$. Now remark that (\ref{String}) -- (\ref{AdaptiveInputs}) may be regarded as the interconnection of the homogeneous stochastic port-Hamiltonian system (\ref{String_PHS}) with the Dirac structure $\mathcal{D}_c$ in the following way
\begin{equation}
    \left\{
    \begin{array}{l}
    f_{\partial,2}(t) = -f_2(t),\\
    e_{\partial,2}(t) = e_2(t),
    \end{array}\right.
    \label{Interconnection}
\end{equation}
where $f_{\partial,2}$ and $e_{\partial,2}$ are the second components of the variables defined in (\ref{FlowEffort}). The interconnection (\ref{Interconnection}) is the same as the one performed in (\ref{InterconnectionDefinition}), which, thanks to Proposition \ref{PropInter}, implies that the controlled and observed system (\ref{String}) -- (\ref{AdaptiveInputs}) may be written as a Dirac structure. In particular, as an interesting feature, the power of the total system (\ref{String}) -- (\ref{AdaptiveInputs}) is zero.

\section{Conclusion \& perspectives}\label{sec:concl}
In this work we introduced and studied the notion of Dirac structure for stochastic port-Hamiltonian systems with multiplicative Gaussian white noise. Taking advantage of the nice geometrical properties of the Stratonovich formalism, the Dirac structure for deterministic infinite-dimensional PHSs as studied in \cite{Kurula} was extended to a stochastic setting. We showed that a newly defined subset of the Cartesian product between extended effort and flow spaces related to a class of SPHSs forms a Dirac structure. As an illustration, we showed that the system composed of a stochastic vibrating string and a mass-spring damper forms a Dirac structure, when interconnected in a power-conserving way. These results should be considered as a first a step towards the development of boundary controllers of SPHSs.

This work opens the way to further research questions and investigations. It would be of great interest to generalize the Dirac structure proposed here for SPHSs by considering various sources of noise entering such as boundary and interconnection noises. Moreover, higher-order stochastic port-Hamiltonian systems will also be considered by the authors in future works.


\section*{Acknowledgments}
This research was conducted with the financial support of F.R.S-FNRS. Anthony Hastir is a FNRS Research Fellow under the grant CR 40010909 and was previously under the grant FC 29535. Francois Lamoline was under the grant FC 08741.
\section*{References}
\bibliographystyle{elsarticle-harv}       
\bibliography{biblio}   

\begin{thebibliography}{22}
\expandafter\ifx\csname natexlab\endcsname\relax\def\natexlab#1{#1}\fi
\providecommand{\url}[1]{\texttt{#1}}
\providecommand{\href}[2]{#2}
\providecommand{\path}[1]{#1}
\providecommand{\DOIprefix}{doi:}
\providecommand{\ArXivprefix}{arXiv:}
\providecommand{\URLprefix}{URL: }
\providecommand{\Pubmedprefix}{pmid:}
\providecommand{\doi}[1]{\href{http://dx.doi.org/#1}{\path{#1}}}
\providecommand{\Pubmed}[1]{\href{pmid:#1}{\path{#1}}}
\providecommand{\bibinfo}[2]{#2}
\ifx\xfnm\relax \def\xfnm[#1]{\unskip,\space#1}\fi
\bibitem[{Caballeria et~al.(2021)Caballeria, Ramirez and Gorrec}]{Caballeria}
\bibinfo{author}{Caballeria, J.}, \bibinfo{author}{Ramirez, H.},
  \bibinfo{author}{Gorrec, Y.L.}, \bibinfo{year}{2021}.
\newblock \bibinfo{title}{An irreversible port-hamiltonian model for a class of
  piezoelectric actuators}.
\newblock \bibinfo{journal}{IFAC-PapersOnLine} \bibinfo{volume}{54},
  \bibinfo{pages}{436--441}.
\newblock \bibinfo{note}{3rd IFAC Conference on Modelling, Identification and
  Control of Nonlinear Systems MICNON 2021}.
\bibitem[{Cordoni et~al.(2019)Cordoni, Persio and
  Muradore}]{cordoni2019stochastic}
\bibinfo{author}{Cordoni, F.}, \bibinfo{author}{Persio, L.D.},
  \bibinfo{author}{Muradore, R.}, \bibinfo{year}{2019}.
\newblock \bibinfo{title}{Stochastic port--hamiltonian systems}.
\newblock \href{http://arxiv.org/abs/1910.01901}{{\tt arXiv:1910.01901}}.
\bibitem[{Cordoni et~al.(2022)Cordoni, Di~Persio and Muradore}]{Cordoni2022}
\bibinfo{author}{Cordoni, F.G.}, \bibinfo{author}{Di~Persio, L.},
  \bibinfo{author}{Muradore, R.}, \bibinfo{year}{2022}.
\newblock \bibinfo{title}{Weak energy shaping for stochastic controlled
  port-hamiltonian systems}.
\newblock \URLprefix \url{https://arxiv.org/abs/2202.08689}.
\bibitem[{Da~Prato and Zabczyk(2014)}]{daprato}
\bibinfo{author}{Da~Prato, G.}, \bibinfo{author}{Zabczyk, J.},
  \bibinfo{year}{2014}.
\newblock \bibinfo{title}{Stochastic Equations in Infinite Dimensions}.
\newblock Encyclopedia of Mathematics and its Applications,
  \bibinfo{publisher}{Cambridge University Press}.
\bibitem[{Duan and Wang(2014)}]{duan2014}
\bibinfo{author}{Duan, J.}, \bibinfo{author}{Wang, W.}, \bibinfo{year}{2014}.
\newblock \bibinfo{title}{Effective Dynamics of Stochastic Partial Differential
  Equations}.
\newblock Elsevier Insights, \bibinfo{publisher}{Elsevier Science}.
\bibitem[{Haddad et~al.(2018)Haddad, Rajpurohit and Jin}]{Haddad}
\bibinfo{author}{Haddad, W.M.}, \bibinfo{author}{Rajpurohit, T.},
  \bibinfo{author}{Jin, X.}, \bibinfo{year}{2018}.
\newblock \bibinfo{title}{Energy-based feedback control for stochastic
  port-controlled hamiltonian systems}.
\newblock \bibinfo{journal}{Automatica} \bibinfo{volume}{97},
  \bibinfo{pages}{134--142}.
\bibitem[{Jacob and Zwart(2012)}]{zwart}
\bibinfo{author}{Jacob, J.}, \bibinfo{author}{Zwart, H.}, \bibinfo{year}{2012}.
\newblock \bibinfo{title}{Linear Port-Hamiltonian Systems on
  Infinite-dimensional Spaces}.
\newblock Number \bibinfo{number}{223} in \bibinfo{series}{Operator Theory:
  Advances and Applications}, \bibinfo{publisher}{Springer},
  \bibinfo{address}{Netherlands}.
\bibitem[{Kurula et~al.(2010)Kurula, Zwart, {van der Schaft} and
  Behrndt}]{Kurula}
\bibinfo{author}{Kurula, M.}, \bibinfo{author}{Zwart, H.},
  \bibinfo{author}{{van der Schaft}, A.}, \bibinfo{author}{Behrndt, J.},
  \bibinfo{year}{2010}.
\newblock \bibinfo{title}{Dirac structures and their composition on hilbert
  spaces}.
\newblock \bibinfo{journal}{Journal of Mathematical Analysis and Applications}
  \bibinfo{volume}{372}, \bibinfo{pages}{402--422}.
\newblock \DOIprefix\doi{https://doi.org/10.1016/j.jmaa.2010.07.004}.
\bibitem[{Lamoline(2019)}]{Lamolinethesis}
\bibinfo{author}{Lamoline, F.}, \bibinfo{year}{2019}.
\newblock \bibinfo{title}{Analysis and {LQG} Control of Infinite-dimensional
  Stochastic Port-{H}amiltonian Systems}.
\newblock Ph.D. thesis. University of Namur.
\bibitem[{Lamoline(2021)}]{LamolineEjc}
\bibinfo{author}{Lamoline, F.}, \bibinfo{year}{2021}.
\newblock \bibinfo{title}{Passivity of boundary controlled and observed
  stochastic port-hamiltonian systems subject to multiplicative and input
  noise}.
\newblock \bibinfo{journal}{European Journal of Control} .
\bibitem[{{Lamoline} and {Winkin}(2017)}]{Lamoline-cdc}
\bibinfo{author}{{Lamoline}, F.}, \bibinfo{author}{{Winkin}, J.J.},
  \bibinfo{year}{2017}.
\newblock \bibinfo{title}{On stochastic port-{H}amiltonian systems with
  boundary control and observation}, in: \bibinfo{booktitle}{2017 IEEE 56th
  Annual Conference on Decision and Control (CDC)}, pp.
  \bibinfo{pages}{2492--2497}.
\bibitem[{{Lamoline} and {Winkin}(2020)}]{LamolineIEEE}
\bibinfo{author}{{Lamoline}, F.}, \bibinfo{author}{{Winkin}, J.J.},
  \bibinfo{year}{2020}.
\newblock \bibinfo{title}{Well-posedness of boundary controlled and observed
  stochastic port-{H}amiltonian systems}.
\newblock \bibinfo{journal}{IEEE Transactions on Automatic Control}
  \bibinfo{volume}{65}, \bibinfo{pages}{4258--4264}.
\bibitem[{{Le Gorrec} et~al.(2005){Le Gorrec}, {Zwart} and
  {Maschke}}]{Legorrec2005}
\bibinfo{author}{{Le Gorrec}, Y.}, \bibinfo{author}{{Zwart}, H.},
  \bibinfo{author}{{Maschke}, B.}, \bibinfo{year}{2005}.
\newblock \bibinfo{title}{Dirac structures and boundary control systems
  associated with skew-symmetric differential operators}.
\newblock \bibinfo{journal}{SIAM Journal on Control and Optimization}
  \bibinfo{volume}{44}, \bibinfo{pages}{1864--1892}.
\bibitem[{Lázaro-Camí and Ortega(2008)}]{Lazaro2008}
\bibinfo{author}{Lázaro-Camí, J.}, \bibinfo{author}{Ortega, J.},
  \bibinfo{year}{2008}.
\newblock \bibinfo{title}{Stochastic {H}amiltonian dynamical systems}.
\newblock \bibinfo{journal}{Reports on Mathematical Physics}
  \bibinfo{volume}{61}, \bibinfo{pages}{65 -- 122}.
\bibitem[{Mora et~al.(2021a)Mora, Gorrec, Ramirez, Yuz and Maschke}]{Mora2021}
\bibinfo{author}{Mora, L.A.}, \bibinfo{author}{Gorrec, Y.L.},
  \bibinfo{author}{Ramirez, H.}, \bibinfo{author}{Yuz, J.},
  \bibinfo{author}{Maschke, B.}, \bibinfo{year}{2021}a.
\newblock \bibinfo{title}{Dissipative port-hamiltonian formulation of maxwell
  viscoelastic fluids}.
\newblock \bibinfo{journal}{IFAC-PapersOnLine} \bibinfo{volume}{54},
  \bibinfo{pages}{430--435}.
\newblock \bibinfo{note}{3rd IFAC Conference on Modelling, Identification and
  Control of Nonlinear Systems MICNON 2021}.
\bibitem[{Mora et~al.(2021b)Mora, Gorrec, Ramírez and Maschke}]{Mora2021Bis}
\bibinfo{author}{Mora, L.A.}, \bibinfo{author}{Gorrec, Y.L.},
  \bibinfo{author}{Ramírez, H.}, \bibinfo{author}{Maschke, B.},
  \bibinfo{year}{2021}b.
\newblock \bibinfo{title}{Irreversible port-hamiltonian modelling of 1d
  compressible fluids}.
\newblock \bibinfo{journal}{IFAC-PapersOnLine} \bibinfo{volume}{54},
  \bibinfo{pages}{64--69}.
\newblock \bibinfo{note}{7th IFAC Workshop on Lagrangian and Hamiltonian
  Methods for Nonlinear Control LHMNC 2021}.
\bibitem[{Ramirez et~al.(2022)Ramirez, Gorrec and Maschke}]{Ramirez2022}
\bibinfo{author}{Ramirez, H.}, \bibinfo{author}{Gorrec, Y.L.},
  \bibinfo{author}{Maschke, B.}, \bibinfo{year}{2022}.
\newblock \bibinfo{title}{Boundary controlled irreversible port-hamiltonian
  systems}.
\newblock \bibinfo{journal}{Chemical Engineering Science}
  \bibinfo{volume}{248}, \bibinfo{pages}{117107}.
\bibitem[{{Rashad Hashem} et~al.(2020){Rashad Hashem}, Califano, {van der
  Schaft} and Stramigioli}]{Hashem20}
\bibinfo{author}{{Rashad Hashem}, R.}, \bibinfo{author}{Califano, F.},
  \bibinfo{author}{{van der Schaft}, A.}, \bibinfo{author}{Stramigioli, S.},
  \bibinfo{year}{2020}.
\newblock \bibinfo{title}{Twenty years of distributed port-hamiltonian systems:
  a literature review}.
\newblock \bibinfo{journal}{IMA journal of mathematical control and
  information} .
\bibitem[{Ruth F.~Curtain(1978)}]{CurtainPritchard}
\bibinfo{author}{Ruth F.~Curtain, A.J.P.}, \bibinfo{year}{1978}.
\newblock \bibinfo{title}{Infinite dimensional linear systems theory}.
\newblock \bibinfo{publisher}{Springer-Verlag}.
\bibitem[{Satoh and Fujimoto(2013)}]{satoh}
\bibinfo{author}{Satoh, S.}, \bibinfo{author}{Fujimoto, K.},
  \bibinfo{year}{2013}.
\newblock \bibinfo{title}{Passivity {B}ased {C}ontrol of {S}tochastic
  {P}ort-{H}amiltonian {S}ystems}.
\newblock \bibinfo{journal}{IEEE Transactions on Automatic Control}
  \bibinfo{volume}{58}, \bibinfo{pages}{1139--1153}.
\bibitem[{van~der Schaft and Maschke(2002)}]{vds02}
\bibinfo{author}{van~der Schaft, A.}, \bibinfo{author}{Maschke, B.},
  \bibinfo{year}{2002}.
\newblock \bibinfo{title}{Hamiltonian formulation of distributed-parameter
  systems with boundary energy flow}.
\newblock \bibinfo{journal}{Journal of Geometry and Physics}
  \bibinfo{volume}{42}, \bibinfo{pages}{166 -- 194}.
\bibitem[{Villegas(2007)}]{Villegas}
\bibinfo{author}{Villegas, J.}, \bibinfo{year}{2007}.
\newblock \bibinfo{title}{A Port-Hamiltonian Approach to distributed parameter
  systems}.
\newblock Ph.D. thesis. University of Twente.

\end{thebibliography}
\end{document}